\documentclass[final]{siamltex}

\usepackage{amsfonts}
\usepackage{amssymb,amsmath}
\usepackage[mathscr]{eucal}
\usepackage{graphicx}
\usepackage{times}

\title{Curvilinearity and Orthogonality}

\author{Vladimir S. Chelyshkov\thanks{Independent Scholar}}

\begin{document}

\maketitle

\begin{abstract}

We introduce sequences of functions orthogonal on a finite interval: proper orthogonal rational functions, orthogonal exponential functions, orthogonal logarithmic functions, and systems of transmuted orthogonal polynomials.
\end{abstract}



\pagestyle{myheadings}
\thispagestyle{plain}\markboth{V.~S.~CHELYSHKOV}{CURVILINEARITY AND ORTHOGONALITY}
\section{Introduction}
General properties of orthogonal rational functions were studied in \cite{C,3C}. Specific examples of orthogonal rational functions on the semi-axis were given in \cite{ 13C, 11C, 12C}. Alternative orthogonal rational functions on a half-line which provide term-by-term increasing rate of decay at infinity were introduced in \cite{4C}.

Orthogonal exponential functions on the semi-axis were presented in \cite{5C}. Alternative orthogonal exponential functions on the semi-axis with term-by-term increasing rate of decay at infinity  were developed in \cite{6C,7C,8C}.

Orthogonal logarithmic functions with weak singularity at one end of an interval were introduced in \cite{SC}.

Standard theory of orthogonal polynomials \cite{9C} includes statement that zeros of the polynomials are distributed in an orthogonality interval. Alternative orthogonal polynomials \cite{7C, 8C} have 
zeros both inside and at one end of an interval. Structured orthogonal polynomials \cite{10C} admit zeros inside and at both ends of an interval. One can also introduce
orthogonal polynomials that have prescribed zeros outside of an orthogonality interval. Being not of specific interest by itself, such altered orthogonal polynomials may originate
orthogonal rational functions and other orthogonal functions on an interval.

In this preliminary notice, we alter Legendre polynomials \cite{2C} and use transforms for independent variable $x$ on the $x$-axis to compose sequences of proper orthogonal rational functions on an interval, as well as sequences of orthogonal exponential functions and orthogonal logarithmic functions on an interval. 

In addition, we introduce systems of transmuted orthogonal polynomials that have complex zeros and may be another origin for above mentioned compositions.

\section{Altered Legendre polynomials} 
Let $\widetilde{P}_n (x)$ be shifted Legendre polynomials.  Making use of $\widetilde{P}_n (x)$, we introduce polynomials 
\begin{equation}
A_{n} (x)=\frac{x+1}{2}\widetilde{P}_{n-1}(x), \quad x \in [-1,1],\quad   n=1,2,…
\label{a1}
\end{equation}
that are orthogonal
\[
\int_0^1 \frac{A_m (x) A_n (x)}{(x+1)^2}{\rm d}x =\frac{\delta_{mn}}{4(2n-1)} 
\]
and have a fixed real zero $x=-1$ outside the interval of orthogonality.
Properties of ${A_{n} (x)}$ can be easily derived from definition (\ref{a1}) and properties of Legendre polynomials.
Altered Legendre polynomials $A_n(x)$ form a basis in $L^2 [0,1]$, and a polynomial $\widetilde{P}_k (x)$ can be approximated by the sequence $A_n (x)$ on the interval.                                                                                                     
In particular, one can find
\begin{equation}
\widetilde{P}_0 (x) =1 \approx 2\ln 2\cdot A_{1}(x) - 6(3\ln 2 -2)A_{2}(x))+10(13\ln2-9)A_{3}(x) - ...\hspace{3pt}.
\label{a3}
\end{equation}
Expansion (\ref{a3}) provides fast convergence for $x\in[0,1]$ and shows good results for extrapolation outside of the interval.

\section{Proper orthogonal rational functions on an interval} Making use of Cayley transform\footnote{The transform permutes elements $\{-1, 0, 1, \infty\}$ in sequence.} \cite{1C}
\begin{equation}
x \mapsto \frac{1-x}{1+x}
\label{b1}
\end{equation}
that maps points $\{ -1,0,1\}$ to points $\{ \infty,1,0\}$ on the $x$-axis, we define proper Legendre  rational functions as a composition
\[
R_n (x):=(-1)^{n-1}A_n \left(\frac{1-x}{1+x}\right),\quad  x \in [0,\infty),\quad  n=1,2, ...\hspace{3pt}.
\]
Functions $R_n(x)$ are orthogonal on the interval $[0,1]$, and
\[
\int_0^1R_m(x)R_n(x){\rm d}x=\frac{\delta_{mn}}{2(2n-1)}.
\]
It may be noted that function 
\[
f(x)=\frac{1-x}{1+x}
\]
is an involution, so function $f^{-1}(x)$ provides the same composition as $f(x)$ does.
\begin{figure}[htbp]
\centerline{\includegraphics[height=36mm, width=60mm]{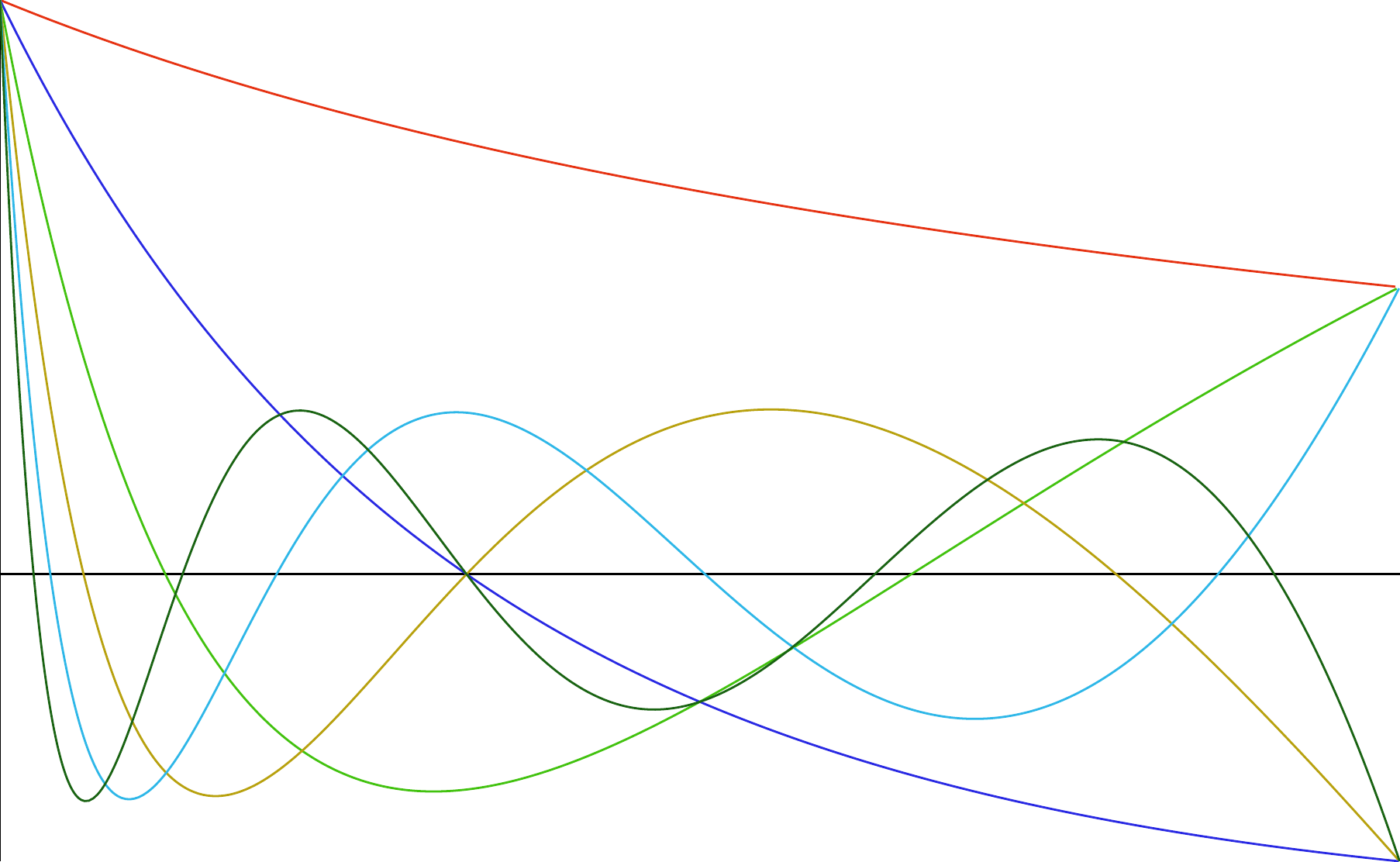}}
\caption{ $R_n(x),\hspace{3pt}n=1-6$.}
\normalsize
\end{figure}
We find
\[
R_1(x)=\frac{1}{x+1},
\]
\[
R_2(x)=\frac{1-3x}{(x+1)^2},
\]
\[
R_3 (x)=\frac{13x^2-10x+1}{(x+1)^3},
\]
\[
R_4(x)=\frac{(1-3x)(21x^2-18x+1)}{(x+1)^4},
\]
\[
R_5 (x)=\frac{321x^4-516x^3+246x^2-36x+1}{(1+x)^5}, 
\]
\[
R_6(x)=\frac{(1-3x)(1-52x+454x^2-948x^3+561x^4)}{(1+x)^6},
\]
and
\[
R_n(0)=1, \quad R_n(1)=\frac{(-1)^{n-1}}{2}, \quad \lim_{x\rightarrow\infty}R_n(x)=0.
\]
From properties of Legendre polynomials it follows that $R_n(x)$ with even indices have a zero at $x=1/3$.

Functions $R_n(x)$ form a basis in $L^2[0,1]$.
In the complex $z$-plain, $R_n(z)$ have one $n$-degree pole at $z=-1$.

\section{Orthogonal exponential functions on an interval} Making use of a transform of exponential type
\begin{equation}
x \mapsto 2^{1-x}-1
\label{b2}
\end{equation}
that maps points $\{-1,0,1\}$  to points  $\{\infty,1,0\}$ on the $x$-axis, we compose a sequence of Legendre exponential functions
\[
E_n(x):=(-1)^{n-1}A_n(2^{1-x}-1),\quad x \in [0, \infty),\quad  n=1,2, ...\hspace{3pt} 
\]
that obey orthogonality relation
\[
\int_0^12^xE_m(x)E_n(x){\rm d}x=\frac{\delta_{mn}}{2\ln2(2n-1)}.
\] 
For convenience, we introduce polynomial functions $F_{n}(t)$ and specify $E_{n}(x)$ as compositions of  $F_{n}(t)$  and an exponential function $t(x)=2^{-x}$, that is,
\[
E_n(x)\equiv F_n(t)_{t=2^{-x}}\quad n=1,2, ...\hspace{3pt}.
\]
\begin{figure}[htbp]
\centerline{\includegraphics[height=36mm, width=60mm]{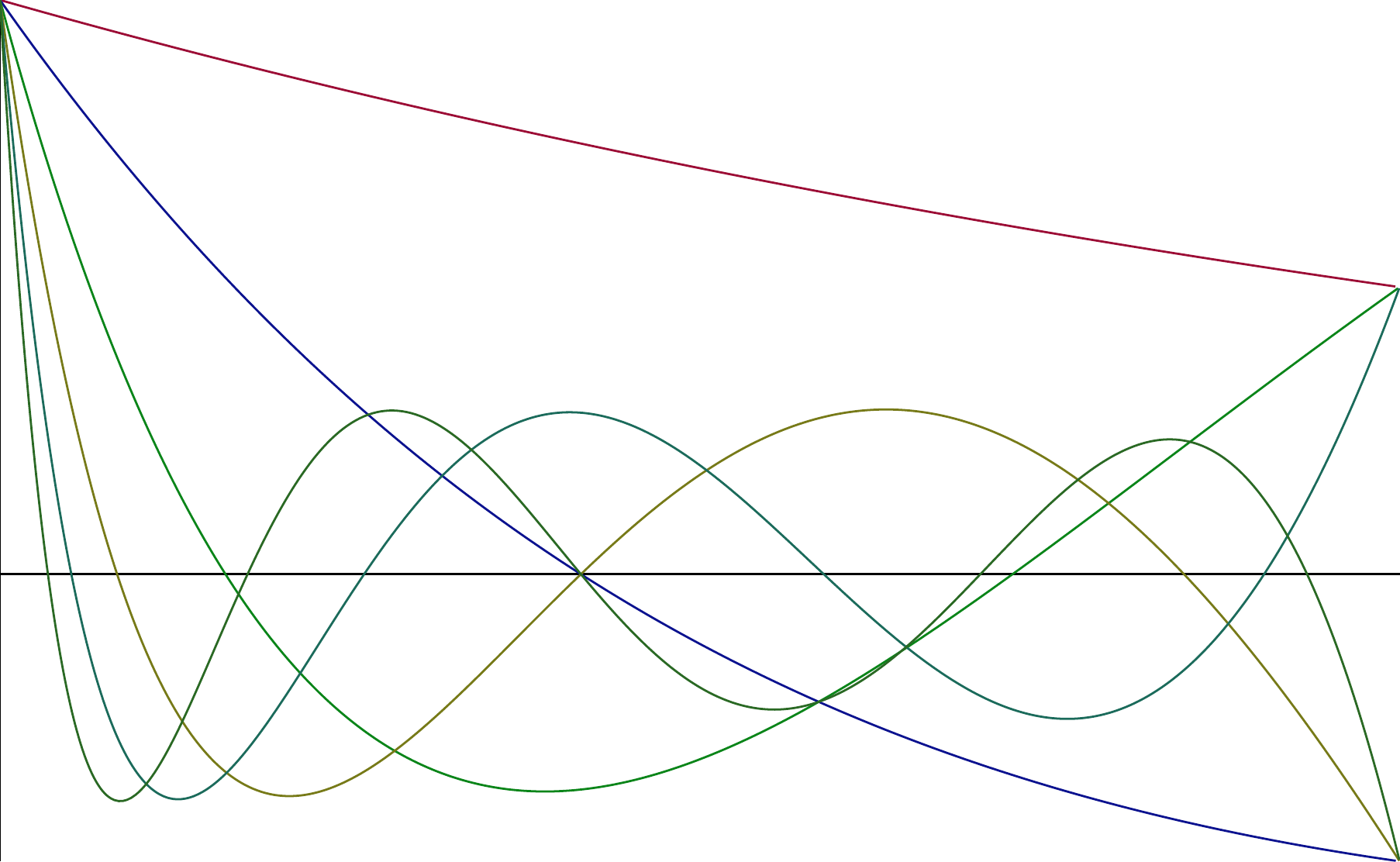}}
\caption{$E_n(x),\hspace{3pt} n=1-6$.}
\normalsize
\end{figure}
We find
\[
F_1 (t)=t,
\] 
\[
F_2 (t)=t(4t-3),
\]
\[
F_3 (t)=t(24t^2-36t+13),
\]
\[
F_4 (t)=t(4t-3)(40t^2-60t+21),
\]
\[
F_5 (t)=t(1120t^4-3360t^3+3720t^2-1800t+321),
\]
\[
F_6 (t)=t(4t-3)(2016t^4-6048t^3+6664t^2-3192t+561),
\]
and
\[
E_n(0)=1, \quad E_n(1)=\frac{(-1)^{n-1}}{2}, \quad \lim_{x\rightarrow\infty}E_n(x)=0.
\]
It appears that $E_n(x)$ with even indices have a zero at
\begin{equation}
x =2 - \ln 3 / \ln 2=0.415037499...\hspace{3pt}.
\label{a4}
\end{equation}

Functions $E_n(x)$ form a basis in $L^2[0,1]$.

\section{Orthogonal logarithmic functions on an interval} A function 
\[
g(x)=2^{1-x}-1.
\]
is not an involution, but it is a near involution on the interval $[0,1]$. Inverse of $g(x)$ is the function  
\[
g^{-1}(x)=1-\log_2(1+x),
\]
and one may find that
\[
\left(\int_0^1(g(x)-g^{-1}(x))^2{\rm d}x\right)^{1/2}=0.002719224...\hspace{3pt}.
\]
Also,
\[
g(x)=g^{-1}(x)\quad {\rm at}\quad x=0, \quad x= 0.456999559..., \quad {\rm and} \quad x=1. 
\]

Making use of a transform of logarithmic type
\begin{equation}
x  \mapsto 1-\log_2(1+x)
\label{b3}
\end{equation}
that maps the points $\{-1, 0, 1\}$ to $\{3, 1, 0\}$ on the $x$-axis, we define Legendre logarithmic functions
\[
 L_n(x):=(-1)^{n-1}A_n(1-\log_2(1+x)),\quad x \in [0, 3],\quad  n=1,2, ...\hspace{3pt}. 
\]
Functions $L_n(x)$ are orthogonal on the interval $[0,1]$,  and
\[
\int_0^1\frac{1}{\frac{1+x}{4}{\rm log}_2^2\frac{1+x}{4}}L_m(x)L_n(x){\rm d}x=\frac{\delta_{mn}}{1/\ln2(2n-1)}.
\]
Also, we introduce polynomial functions $G_{n}(t)$ and identify $L_{n}(x)$ as compositions of $G_n(t)$ and a logarithmic function $t(x)=\log_2(1+x)$,
\[
L_n(x)\equiv G_n(t)_{t=\log_2(1+x)}\quad n=1,2, ...\hspace{3pt}.
\]
We find
\[
G_n(t) = (-1)^{n-1}(1-\frac{1}{2}t)\widetilde{P}_{n-1}(t),
\]
and
\[
L_n(0)=1, \quad L_n(1)=\frac{(-1)^{n-1}}{2}, \quad L_n(3)=0.
\]
It appears that $L_n(x)$ with even indices have a zero at
\begin{equation}
x=\sqrt{2}-1=0.414213562...\hspace{3pt}.
\label{a5}
\end{equation}
Numbers in (\ref{a4}) and (\ref{a5}) are ``almost the same''. Actually,
\begin{figure}[htbp]
\centerline{\includegraphics[height=36mm, width=60mm]{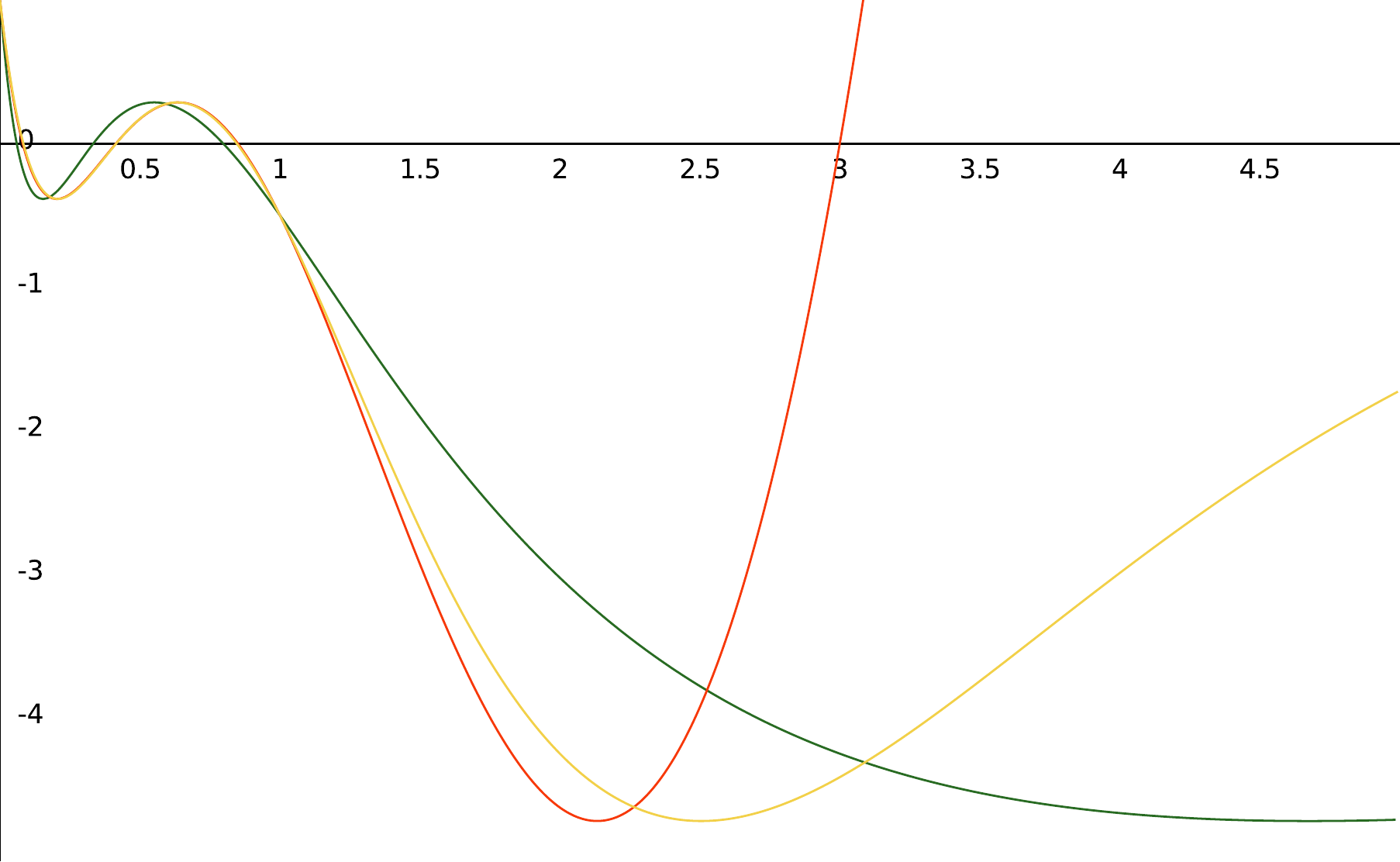}}
\caption{ $R_4(x)$ (green), $E_4(x)$ (yellow), $L_4(x)$ (red).}
\normalsize
\end{figure}
the graphs of the functions $E_n(x)$ and $L_n(x)$ are very close on the interval of orthogonality, but they differ outside of the interval.

Functions $L_n(x)$ form a basis in $L^2[0,1]$.

\section{Conjecture}Calculations show that the following statement may take place. Let $r_{nk}$, $e_{nk}$, and $l_{nk}$ be abscissas of local extrema for functions $R_{n}(x)$, $E_{n}(x)$, and $L_{n}(x)$. Then $R_{n}(r_{nk})=E_{n}(e_{nk})=L_{n}(l_{nk})$. 

\section{Conclusion} Altered Legendre polynomials ${A}_{n}(x)$ originate three sequences of orthogonal functions, $R_{n}(x)$, $E_n(x)$, and $L_n(x)$.

The sequences are bases  for approximation in $L^2[0,1]$ space, and they may provide spectral convergence properties.  One may notice that a linear function $s(x)=ax+b$ cannot be represented exactly by a finite linear combination of any of the three sequences, and the bases do not provide options for linear approximation. 

Similar compositions for altered alternative orthogonal polynomials can be considered.

Also, one may extend definition of $A_n(x)$ to a three-parametric family of functions
\[
_{\gamma}{A}^{(\alpha,\beta)}_{n} (x)=((x+1)/2)^\gamma\widetilde{P}_{n-1}^{(\alpha,\beta)}(x), \quad \gamma>0,\quad n>0,
\]
where $\widetilde{P}_{n-1}^{(\alpha,\beta)}(x)$ are shifted Jacobi  polynomials.

\section{Appendix 1. Orthogonal exponential functions on an interval, the general case} 
 We introduce a parametric transform of exponential type
\begin{equation}
x \mapsto \frac{b^{1-x}-1}{b-1}, \quad b  >1,
\label{nu1}
\end{equation}
and define a single-parameter $b$-sequense of alterd Legendre polynomials
\begin{equation}
A_n(b; x):=\alpha(x+a)\widetilde{P}_{n-1} (x), \quad x \in [-a,1],\quad   n=1,2,…\hspace{2pt},
\label{nu0}
\end{equation}
with
\[
\alpha =(b-1)/b,\quad \quad a=1/(b-1).
\]
Transform (\ref{nu1}) maps points $\{-a,0,1\}$  to points  $\{\infty,1,0\}$ on the $x$-axis, and we compose a  $b$-sequence of Legendre exponential functions
\[
E_n(b; x):=(-1)^{n-1}A_n \left(b;  \frac{b^{1-x}-1}{b-1}\right), \quad x \in [0, \infty),\quad   n=1,2,…
\]
that obey orthogonality relations
\[
\int_0^1b^xE_m(b; x)E_n(b; x){\rm d}x=\frac{(b-1)\delta_{mn}}{b\ln(b)(2n-1)}.
\]
We find
\[
E_1(b; x)=b^{-x},
\]
\[
E_2(b; x)=\frac{2bb^{-2x}-(b+1)b^{-x}}{b-1},
\]
\[
E_3(b; x)=\frac{6b^2b^{-3x}-6b(b+1)b^{-2x}+(b^2+4b+1)b^{-x}}{(b-1)^2},
\]
\[
E_4(b; x)=\frac{20b^3b^{-4x}-30b^2(b+1)b^{-3x}+12b(b^2+3b+1)b^{-2x}}{(b-1)^3}
\]
\[
-\frac{(b^3+9b^2+9b+1)b^{-x}}{(b-1)^3},
\]
\[
E_5(b; x)=\frac{70b^4b^{-5x}-140b^3(b+1)b^{-4x}+30b^2(3b^2+8b+3)b^{-3x}}{(b-1)^4}
\]
\[
+\frac{-20b(b^3+6b^2+6b+1)b^{-2x}+(b^4+16b^3+36b^2+16b+1)b^{-x}}{(b-1)^4},
\]
\[
E_6(b; x)=\frac{252b^5b^{-6x}-630b^4(b+1)b^{-5x}+280b^3(2b^2+5b+2))b^{-4x}}{(b-1)^5}
\]
\[
+\frac{-210b^2(b^3+5b^2+5b+1)b^{-3x}+30b(b^4+10b^3+20b^2+10b+1)b^{-2x}}{(b-1)^5}
\]
\[
-\frac{(b^5+25b^4+100b^3+100b^2+25b+1)b^{-x}}{(b-1)^5},
\]
and
\[
E_n(b; 0)=1, \quad E_n(b; 1)=\frac{(-1)^{n-1}}{b}, \quad \lim_{x\rightarrow\infty}E_n(b; x)=0.
\]
Also, for even $n$ and $x=\log_{b}\frac{2b}{b+1}$, $E_n(b; x)=0$,\quad for $b=e,\hspace{3pt} x= 0.379885$... .

Functions  $E_n(b; x)$ form a basis in $L_2[0,1]$.

\section{Appendix 2. Orthogonal logarithmic functions on an interval, the general case}

Let
\[
h(b; x) = \frac{b^{1-x}-1}{b-1}, \quad b >1,
\]
be a function of variable $x$ and parameter $b$. Inverse to $h(b; x)$ is the function
\[
h^{-1}(b; x)=1-{\ell}_b(x),
\]
\[
{\ell}_b(x)=\log_b(1+(b-1)x).
\]
Making use of a parametric transform of logarithmic type
\begin{equation}
x  \mapsto h^{-1}(b; x)
\label{9.1}
\end{equation}
that maps points $\{-a, 0, 1\}$ to $\{\beta, 1, 0\}$ with 
$\beta=(b^{b/(b-1)}-1)/(b-1)$ 
we define\\ $b$~-~parametric Legendre logarithmic functions
\[
L_n(b; x):=(-1)^{n-1}A_n(b; 1-{\ell}_b(x)),\quad x \in [0, \beta],\quad  n=1,2, ...\hspace{3pt}. 
\]
Functions $L_n(b; x)$ are orthogonal on the interval $[0,1]$, and 
\[
\int_0^1\frac{1}{(1/(b-1)+x)(1-\ell_b(x)/b)^2}L_m(b; x)L_n(b; x)d x =\frac{(b-1)\delta_{mn}}{1/\ln b(2n-1)}.
\]
We find
\[
L_1(b; x) = 1-(b-1)/b \ell_{b}(x),
\]
\[
L_2(b; x) = (1-(b-1)/b\ell_{b}(x))(1-2\ell_b(x)),
\]
and, in general,
\[
L_n(b; x):=(-1)^{n-1}(1-(b-1)/b\ell_{b}(x))\widetilde{P}_{n-1}(1-\ell_{b}(x)),
\]
\[
L_n(b; 0)=1, \quad L_n(b; 1)=\frac{(-1)^{n-1}}{b}, \quad  L_n(b; \beta)=0.
\]
For even $n$, $L_n(b; x)=0$ at $x=(\sqrt{b}-1)/(b-1)$. Also, functions $L_n(b; x)$ are complex-valued for $x< -1/(b-1)$.

Functions $h(b; x)$ and $h^{-1}(b; x)$ are close to involution  for $b=2$ and $b=e$  on the orthogonality interval.

Functions $L_n(b; x)$ form a basis in $L_2[0, 1]$.

\section{Appendix 3. Proper orthogonal rational functions on an interval, a single-parameter case}  We introduce a parametric linear fractional  transformation
\begin{equation}
x \mapsto \frac{1-x}{1+cx}, \quad c  >0,
\label{nunu1}
\end{equation}
and define a single-parameter $c$-sequense of alterd Legendre polynomials (cf.(\ref{nu0}))
\[
A_n(c; x):=\alpha(x+a)\widetilde{P}_{n-1} (x), \quad x \in [-a,1],\quad   n=1,2,…\hspace{2pt},
\]
with
\[
\alpha =c/(c+1),\quad \quad a=1/c.
\]
Transformation (\ref{nunu1}) maps points $\{-a,0,1\}$  to points  $\{\infty,1,0\}$ on the $x$-axis, and we compose a  $c$-sequence of Legendre rational functions
\[
R_n(c; x):=(-1)^{n-1}A_n \left(c; \frac{1-x}{1+cx}\right), \quad x \in [0, \infty),\quad   n=1,2,…
\]
that obey orthogonality relations
\[
\int_0^1R_m(c; x)R_n(c; x){\rm d}x=\frac{\delta_{mn}}{(1+c)(2n-1)}.
\]
We find
\[
R_1(c; x)=\frac{1}{1+cx},
\]
\[
R_2(c; x)=-\frac{(c+2)x-1}{(1+cx)^2},
\]
\[
R_3(c; x)=\frac{(c^2+6c+6)x^2-2(2c+3)x+1}{(1+cx)^3},
\]
\[
R_4(c; x)=-\frac{(c^3+12c^2+30c+20)x^3+(-9c^2-36c-30)x^2+(9c+12)x-1}{(1+cx)^4}
\]
\[
\quad\equiv R_2(c; x)\frac{(c^2+10c+10)x^2-(8c+10)x+1}{(1+cx)^2},
\]
and
\[
R_n(c; 0)=1, \quad R_n(c; 1)=\frac{(-1)^{n-1}}{1+c}, \quad \lim_{x\rightarrow\infty}R_n(c; x)=0.
\]
From properties of Legendre polynomials it follows that $R_n(c; x)$ with even indices have a zero at $x=1/(c+2)$.

Functions  $R_n(c; x)$ form a basis in $L^2[0,1]$.
In the complex $z$-plain, $R_n(c; z)$ have one $n$-degree pole at $z=-1/c$.

\section{Appendix 4. Transmuted orthogonal polynomials} In the previous sections, we used transforms (\ref{b1}), (\ref{b2}), (\ref{b3}), (\ref{nu1}), (\ref{9.1}) on the $x$-axis for composing five sequences of non-polynomial orthogonal functions on the interval $[0,1]$. The transforms represent $x$-mapping of the interval on itself,
\[
[0,1] \mapsto [1,0],
\]
which also can be defined by any term of a sequence of polynomial mappings 
\begin{equation}
x \mapsto 1+\sum_{j=1}^{\nu} (-x)^{j},\quad\
\nu=2m-1,\quad m=1,2,...\hspace{3pt}.
\label{a6}
\end{equation} 
The terms of sequence (\ref{a6}) can be applied for transmuting orthogonal polynomials. In particular, one can transmute shifted Legendre polynomials  $\widetilde{P}_n (x)$ to a system of orthogonal polynomials 
\begin{equation}
\mathcal{P}_{\nu0}(x):=1, \quad \mathcal{P}_{{\nu}n}(x):=\widetilde{P}_n\left (1+\sum_{j=1}^{\nu} (-x)^{j}\right), \quad n>0 
\label{a7}
\end{equation}
of degree ${\nu}n$ that obey orthogonality relations
\[
\int_0^1w_{\nu}(x)\mathcal{P}_{{\nu}k}(x)\mathcal{P}_{{\nu}l}(x){\rm d}x=\frac{\delta_{kl}}{2k+1},
\]
where 
\[
w_{\nu}(x)=\sum_{j=1}^{\nu} j(-x)^{j-1}>0,\quad w_{\nu}(0)=1,\quad w_{\nu}(1)=m
\]
is a polynomial weight function for a chosen $\nu$. 
\begin{figure}[htbp]
\centerline{\includegraphics[height=36mm, width=60mm]{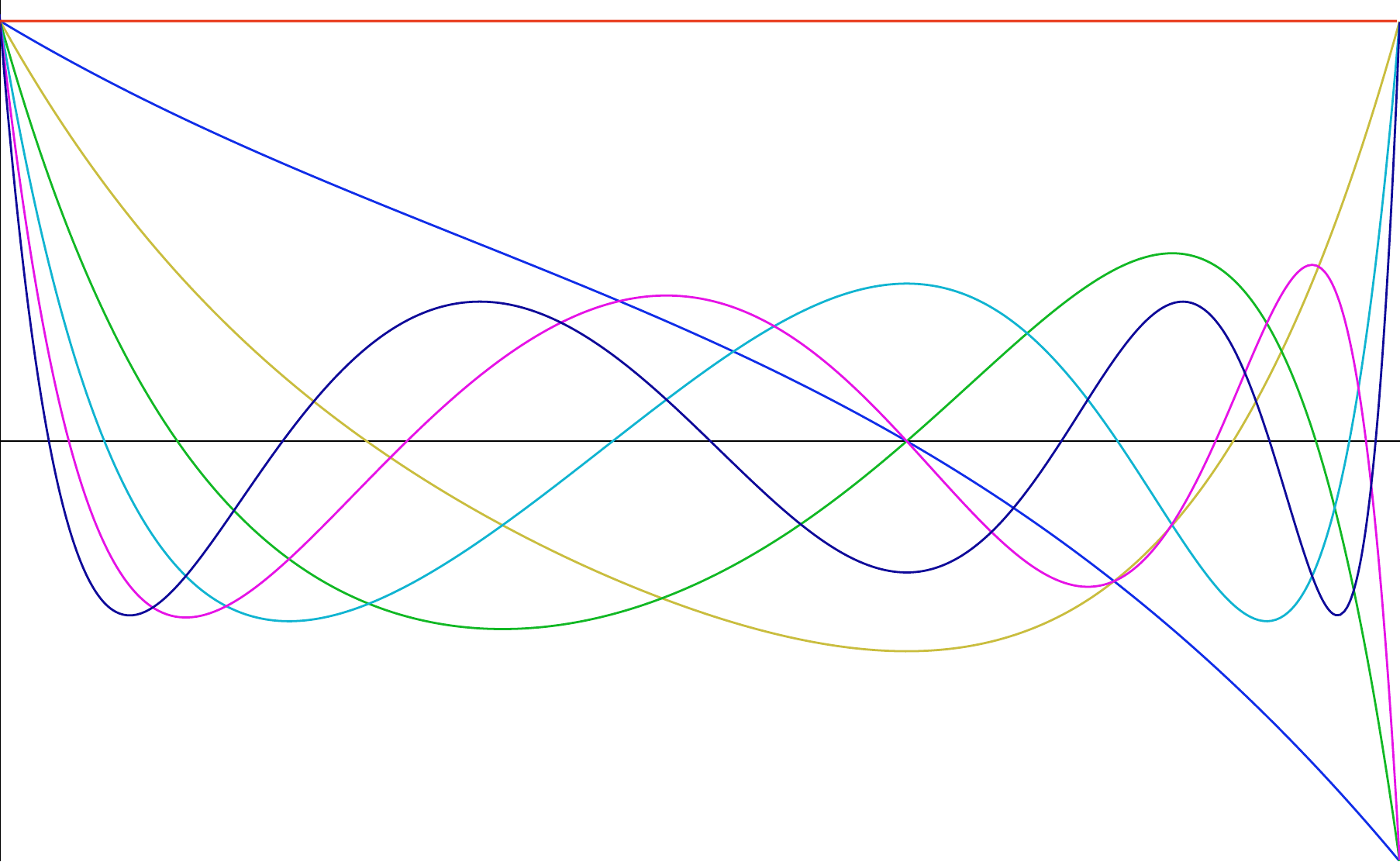}}
\caption{$\mathcal{P}_{3n}(x),\hspace{3pt} n=0-6$.}
\normalsize
\end{figure}

For ${\nu}=3$ polynomial composition (\ref{a7}) results in sequence
\[
\mathcal{P}_{30}(x)=1,
\]
\[
\mathcal{P}_{31}(x)=-2x^3+2x^2-2x+1,
\]
\[
\mathcal{P}_{32}(x)=6x^6-12x^5+18x^4-18x^3+12x^2-6x+1,
\]
\[
\mathcal{P}_{33}(x)=-20x^9+60x^8-120x^7+170x^6-180x^5+150x^4-92x^3+42x^2-12x+1,
\]
\[
\mathcal{P}_{34}= (x)=70x^{12}-280x^{11}+700x^{10}-1260x^9+1750x^8-1960x^7+1770x^6
\]
\[
-1300x^5+760x^4-340x^3+110x^2-20x+1,
\]
\[
\mathcal{P}_{35}(x)=252x^{15}-1260x^{14}+3780x^{13}-8190x^{12}+13860x^{11}-19152x^{10}+21980x^9
\]
\[
-21210x^8+17220x^7-11690x^6+6552x^5-2940x^4+1010x^3-240x^2+30x+1, 
\]
\[
\mathcal{P}_{36}(x)=
924x^{18}-5544x^{17}+19404x^{16}-48972x^{15}+97020x^{14}-158004x^{13}+216594x^{12}
\]
\[
-253764x^{11}+256032x^{10}-223020x^9+167454x^8-107604x^7+58464x^6
\]
\[
-26292x^5+9450x^4-2562x^3+462x^2-42x+1.
\]

     From properties of Legendre polynomials it follows that $\mathcal{P}_{\nu n}(x)$ with $\nu=3$ and any odd $n$ have a zero at $x=0.647798871...$\hspace{2pt}. Calculations show that an increase in 
$\nu$ leads to an increase in this value of $x$ in the interval $[0,1]$, and the zeros cluster towards $x=1$.

Polynomials $\mathcal{P}_{{\nu}n}(x)$ have $n$ real zeros and $(\nu-1)n$ complex conjugate zeros. We find that transmutation (\ref{a6}) results in redistribution of real zeros of the standard orthogonal polynomials and parameter $\nu$ can be selected to increase density of the zeros near $x=1$.  

Transmuted orthogonal polynomials are bases in $L_2[0,1]$. The polynomials originate compositions that are similar to those  introduced in sections  2 - 5, 8 - 10.

\end{document}